\documentclass[11pt]{article}
\usepackage[a4paper,margin=2.5truecm]{geometry}
\usepackage{amsmath,amsthm,amsfonts,amssymb,hyperref,graphicx,t1enc}
\newtheorem{definition}{Definition}
\newtheorem{corollary}{Corollary}
\newtheorem{lemma}{Lemma}
\newtheorem{theorem}{Theorem}
\newcommand{\N}{\mathbb N}
\newcommand{\be}{\beta}
\newcommand{\wt}{\widetilde}

\title{The set \mbox{$\{2^{\textstyle p} \cdot 3^{\textstyle q}:~(p,q \in \N \setminus \{0,1\})~\wedge$}
\mbox{$(\forall x_2,\ldots,x_p \in \N \setminus \{0,1\}~\exists y_2,\ldots,y_p \in \{2,\ldots,q\}$}
\mbox{$\forall i,j,k \in \{2,\ldots,p\}~(x_i^{\textstyle x_j}=x_k \Rightarrow y_i^{\textstyle y_j}=y_k))\}$}
\mbox{is not recursively enumerable}}
\author{Apoloniusz Tyszka}
\begin{document}
\date{}
\begin{sloppypar}
\maketitle
\begin{abstract}
We prove that the set
\[
W=\{\frac{(p+q)(p+q+1)}{2}+q: (p,q \in \N) ~\wedge 
\]
\[
\forall (x_0,\ldots,x_p) \in {\N}^{p+1} ~\exists (y_0,\ldots,y_p) \in \{0,\ldots,q\}^{p+1}
\]
\[
((\forall j,k \in \{0,\ldots,p\} ~(x_j+1=x_k \Rightarrow y_j+1=y_k))~\wedge
\]
\[
(\forall i,j,k \in \{0,\ldots,p\} ~(x_i \cdot x_j=x_k \Rightarrow y_i \cdot y_j=y_k)))\}
\]
is co-recursively enumerable and not recursively enumerable.
Let \mbox{$\be:{\N}^3 \to \N$} denote G\"odel's $\be$ function.
For \mbox{$x_1,x_2,x_3 \in \N$},
\mbox{$\be(x_1,x_2,x_3)$} equals the remainder after integer division of $x_1$ by \mbox{$1+(x_3+1) \cdot x_2$}.
We prove that the set $W$ consists of all \mbox{$n \in \N$} such that
\[
\forall u,v \in \N ~\exists a,b,p,q \in \N ~((2n=(p+q)(p+q+1)+2q) ~\wedge ~\forall i,j,k \in \{0,\ldots,p\}
\]
\[
((\be(a,b,i) \leqslant q) ~\wedge ~(\be(u,v,j)+1=\be(u,v,k) \Rightarrow \be(a,b,j)+1=\be(a,b,k)) ~\wedge
\]
\[
(\be(u,v,i) \cdot \be(u,v,j)=\be(u,v,k) \Rightarrow \be(a,b,i) \cdot \be(a,b,j)=\be(a,b,k))))
\]
We express the above formula in Peano arithmetic.
\end{abstract}
\mbox{limit-computable} function, Peano arithmetic, recursively enumerable set.
\vskip 1.0truecm
\section{\normalsize Introduction}
Let \mbox{$F(x,n)$} denote the formula
\[
\exists ab ~\forall i \leqslant n ~\exists swpq ~\forall jv ~\exists eg ~\{(s+w)^2+3w+s=2i ~\wedge ~\langle[j=w ~\vee ~v=q]
\]
\[
\vee~[j=3i ~\wedge ~v=p+q] ~\vee ~[j=s ~\wedge ~(v=p ~\vee ~(i=n ~\wedge ~v=q+x))]
\]
\[
\vee~[j=3i+1 ~\wedge ~v=pq] ~\Rightarrow ~a=v+e+ejb ~\wedge ~v+g=jb\rangle\}
\]
from \cite[p.~336]{Jones}. From the results of \cite{Jones}, it follows that the set \mbox{$J=\{n \in \N: \neg F(n,n)\}$}
is co-recursively enumerable and not recursively enumerable. We prove that the set
\[
\{2^{\textstyle p} \cdot 3^{\textstyle q}: (p,q \in \N \setminus \{0,1\}) ~\wedge
\]
\[
(\forall x_2,\ldots,x_p \in \N \setminus \{0,1\} ~\exists y_2,\ldots,y_p \in \{2,\ldots,q\}
~\forall i,j,k \in \{2,\ldots,p\} ~(x_i^{\textstyle x_j}=x_k \Rightarrow y_i^{\textstyle y_j}=y_k))\}
\]
is co-recursively enumerable and not recursively enumerable.
For the author, it is the simplest known example of a \mbox{non-recursively} enumerable set \mbox{$X \subseteq \N$}.
\section{\normalsize Limit-computable functions}
Semi-algorithms differ from algorithms, as they may not terminate.
\begin{definition}\label{defi0} (cf.~\cite[pp.~233--235]{Soare}).
A computation in the limit of a function \mbox{$f:\N \to \N$} is a \mbox{semi-algorithm}
which takes as input a \mbox{non-negative} integer $n$ and for every \mbox{$m \in \N$} prints a \mbox{non-negative}
integer \mbox{$\xi(n,m)$} such that \mbox{$\lim \limits_{m \to \infty} \xi(n,m)=f(n)$}.
\end{definition}
\par
By Definition~\ref{defi0}, a function \mbox{$f:\N \to \N$} is computable in the limit
when there exists an infinite computation which takes as input a non-negative integer $n$
and prints a non-negative integer on each iteration and prints \mbox{$f(n)$} on each
sufficiently high iteration.
\section{\normalsize A limit-computable function \mbox{$f:\N \to \N$} which eventually dominates
every computable function \mbox{$g:\N \to \N$}}
For \mbox{$n \in \N$}, let
\[
E_n=\{1=x_k,~x_i+x_j=x_k,~x_i \cdot x_j=x_k:~i,j,k \in \{0,\ldots,n\}\}
\]
\par
For \mbox{$n \in \N$}, \mbox{$f(n)$} denotes the smallest \mbox{$b \in \N$} such that
if a system of equations \mbox{$S \subseteq E_n$} has a solution in \mbox{${\N}^{n+1}$},
then $S$ has a solution in \mbox{$\{0,\ldots,b\}^{n+1}$}.
The function \mbox{$f:\N \to \N$} is computable in the limit and
eventually dominates every computable function \mbox{$g:\N \to \N$}, see \cite{Tyszka1}.
The term {\em "dominated"} in the title of \cite{Tyszka1} means {\em "eventually dominated"}.
\begin{theorem}\label{theo2} (\cite{Tyszka1}).
Flowchart 1 shows a \mbox{semi-algorithm} which computes \mbox{$f(n)$} in the limit.
\end{theorem}
\begin{center}
\includegraphics[scale=0.9167]{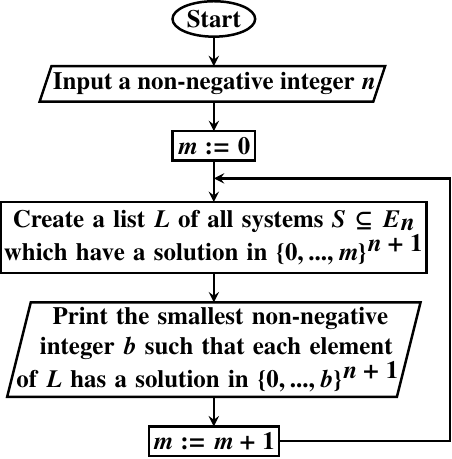}
\end{center}
\vskip 0.01truecm
\centerline{{\bf Flowchart 1}}
\vskip 0.1truecm
\centerline{A \mbox{semi-algorithm} which computes \mbox{$f(n)$} in the limit}
\begin{definition}\label{defi1}
An approximation of a tuple \mbox{$(x_0,\ldots,x_n) \in {\N}^{n+1}$} is a tuple
\mbox{$(y_0,\ldots,y_n) \in {\N}^{n+1}$} such that
\[
(\forall k \in \{0,\ldots,n\} ~(1=x_k \Rightarrow 1=y_k))~\wedge
\]
\[
(\forall i,j,k \in \{0,\ldots,n\} ~(x_i+x_j=x_k \Rightarrow y_i+y_j=y_k))~\wedge
\]
\[
(\forall i,j,k \in \{0,\ldots,n\} ~(x_i \cdot x_j=x_k \Rightarrow y_i \cdot y_j=y_k))
\]
\end{definition}
\newpage
Flowchart 2 shows a simpler semi-algorithm which computes \mbox{$f(n)$} in the limit, see Theorem~\ref{theo3}.
\begin{center}
\includegraphics[scale=0.9167]{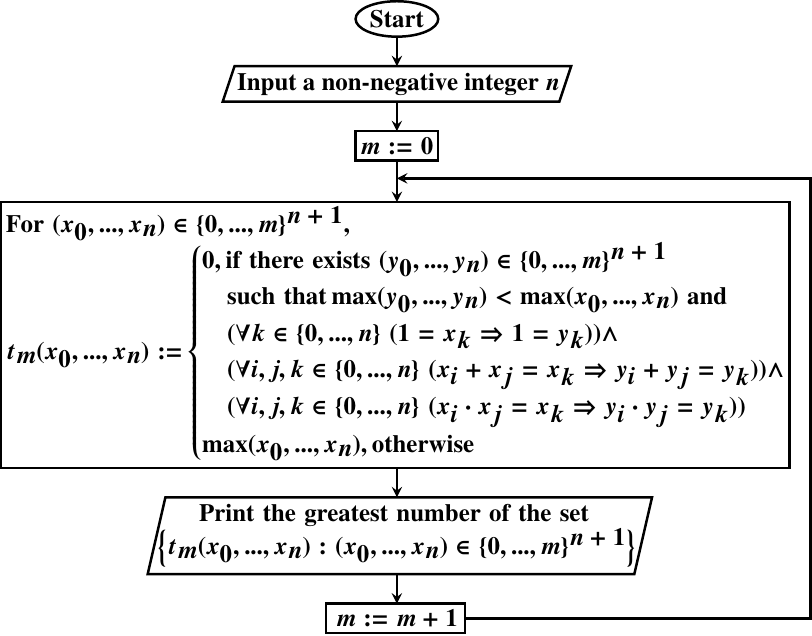}
\end{center}
\vskip 0.01truecm
\centerline{{\bf Flowchart 2}}
\vskip 0.1truecm
\centerline{A simpler \mbox{semi-algorithm} which computes \mbox{$f(n)$} in the limit}
\begin{lemma}\label{lemm1}
For every \mbox{$n,m \in \N$}, the number printed by Flowchart 2 does not exceed the number printed by Flowchart 1.
\end{lemma}
\begin{proof}
For every \mbox{$(a_0,\ldots,a_n) \in \{0,\ldots,m\}^{n+1}$},
\[
E_n \supseteq \{1=x_k:~(k \in \{0,\ldots,n\}) \wedge (1=a_k)\}~\cup
\]
\[
\{x_i+x_j=x_k:~(i,j,k \in \{0,\ldots,n\}) \wedge (a_i+a_j=a_k)\}~\cup
\]
\[
\{x_i \cdot x_j=x_k:~(i,j,k \in \{0,\ldots,n\}) \wedge (a_i \cdot a_j=a_k)\}
\]
\end{proof}
\begin{lemma}\label{lemm2}
For every \mbox{$n,m \in \N$}, the number printed by Flowchart 1 does not exceed the number printed by Flowchart 2.
\end{lemma}
\begin{proof}
Let \mbox{$n,m \in \N$}. For every system of equations \mbox{$S \subseteq E_n$},
if \mbox{$(a_0,\ldots,a_n) \in \{0,\ldots,m\}^{n+1}$} and \mbox{$(a_0,\ldots,a_n)$} solves $S$,
then \mbox{$(a_0,\ldots,a_n)$} solves the following system of equations:
\[
\{1=x_k:~(k \in \{0,\ldots,n\}) \wedge (1=a_k)\}~\cup
\]
\[
\{x_i+x_j=x_k:~(i,j,k \in \{0,\ldots,n\}) \wedge (a_i+a_j=a_k)\}~\cup
\]
\[
\{x_i \cdot x_j=x_k:~(i,j,k \in \{0,\ldots,n\}) \wedge (a_i \cdot a_j=a_k)\}
\]
\end{proof}
\newpage
\begin{theorem}\label{theo3}
For every \mbox{$n,m \in \N$}, Flowcharts 1 and 2 print the same number.
\end{theorem}
\begin{proof}
It follows from Lemmas~\ref{lemm1} and~\ref{lemm2}.
\end{proof}
\begin{corollary}\label{coro1}
For every \mbox{$n,m \in \N$}, Flowcharts 1 and 2 print the smallest \mbox{$b \in \{0,\ldots,m\}$} such that every tuple
\mbox{$(x_0,\ldots,x_n) \in \{0,\ldots,m\}^{n+1}$} has an approximation in \mbox{$\{0,\ldots,b\}^{n+1}$}.
\end{corollary}
\begin{theorem}\label{theo4}
For every \mbox{$n \in \N$}, \mbox{$f(n)$} is the smallest \mbox{$b \in \N$} such that every tuple
\mbox{$(x_0,\ldots,x_n) \in {\N}^{n+1}$} has an approximation in \mbox{$\{0,\ldots,b\}^{n+1}$}.
\end{theorem}
\begin{proof}
It follows from Theorem~\ref{theo2} and Corollary~\ref{coro1}.
\end{proof}
\section{\normalsize Two undecidable decision problems which are related to the function $f$}
\begin{theorem}\label{theo5}
No algorithm takes as input \mbox{non-negative integers} $n$ and $m$ and decides whether or not
\[
\forall (x_0,\ldots,x_n) \in {\N}^{n+1} ~\exists (y_0,\ldots,y_n) \in \{0,\ldots,m\}^{n+1}
\]
\[
((\forall k \in \{0,\ldots,n\} ~(1=x_k \Rightarrow 1=y_k))~\wedge
\]
\[
(\forall i,j,k \in \{0,\ldots,n\} ~(x_i+x_j=x_k \Rightarrow y_i+y_j=y_k))~\wedge
\]
\[
(\forall i,j,k \in \{0,\ldots,n\} ~(x_i \cdot x_j=x_k \Rightarrow y_i \cdot y_j=y_k)))
\]
\end{theorem}
\begin{proof}
Since the function $f$ is not computable, it follows from Theorem~\ref{theo4}.
\end{proof}
\begin{lemma}\label{Cantor} (\cite{Cantor}).
The function
\[
{\N}^2 \ni (p,q) \to \frac{(p+q)(p+q+1)}{2}+q \in \N
\]
is bijective.
\end{lemma}
\begin{theorem}\label{theo6}
No algorithm takes as input a non-negative integer $n$ and decides whether or not
\[
\exists p,q \in \N ~((2n=(p+q)(p+q+1)+2q))~\wedge
\]
\[
\forall (x_0,\ldots,x_p) \in {\N}^{p+1} ~\exists (y_0,\ldots,y_p) \in \{0,\ldots,q\}^{p+1}
\]
\[
((\forall k \in \{0,\ldots,p\} ~(1=x_k \Rightarrow 1=y_k))~\wedge
\]
\[
(\forall i,j,k \in \{0,\ldots,p\} ~(x_i+x_j=x_k \Rightarrow y_i+y_j=y_k))~\wedge
\]
\[
(\forall i,j,k \in \{0,\ldots,p\} ~(x_i \cdot x_j=x_k \Rightarrow y_i \cdot y_j=y_k))))
\]
\end{theorem}
\begin{proof}
It follows from Theorem~\ref{theo5} and Lemma~\ref{Cantor}.
\end{proof}
\section{\normalsize A simple example of a non-recursively enumerable set \mbox{$T \subseteq \N$}, where\\
$T=\{\frac{(p+q)(p+q+1)}{2}+q: (p,q \in \N) \wedge \Psi(p,q)\}$
and the formula $(p,q \in \N) \wedge \Psi(p,q)$ concerns $({\N}^{p+1},\{0,\ldots,q\}^{p+1})$}
Let
\[
T=\{n \in \N:~\exists p,q \in \N ~((2n=(p+q)(p+q+1)+2q)~\wedge
\]
\[
\forall (x_0,\ldots,x_p) \in {\N}^{p+1} ~\exists (y_0,\ldots,y_p) \in \{0,\ldots,q\}^{p+1}
\]
\[
((\forall k \in \{0,\ldots,p\} ~(1=x_k \Rightarrow 1=y_k))~\wedge
\]
\[
(\forall i,j,k \in \{0,\ldots,p\} ~(x_i+x_j=x_k \Rightarrow y_i+y_j=y_k))~\wedge
\]
\[
(\forall i,j,k \in \{0,\ldots,p\} ~(x_i \cdot x_j=x_k \Rightarrow y_i \cdot y_j=y_k))))\}
\]
\begin{theorem}\label{theo7}
The set \mbox{$\N \setminus T$} is recursively enumerable.
\end{theorem}
\begin{proof}
Flowchart 3 shows a \mbox{semi-algorithm} which takes as input \mbox{$n \in \N$}
and terminates if and only if \mbox{$n \in \N \setminus T$}.
\begin{center}
\includegraphics[scale=0.9167]{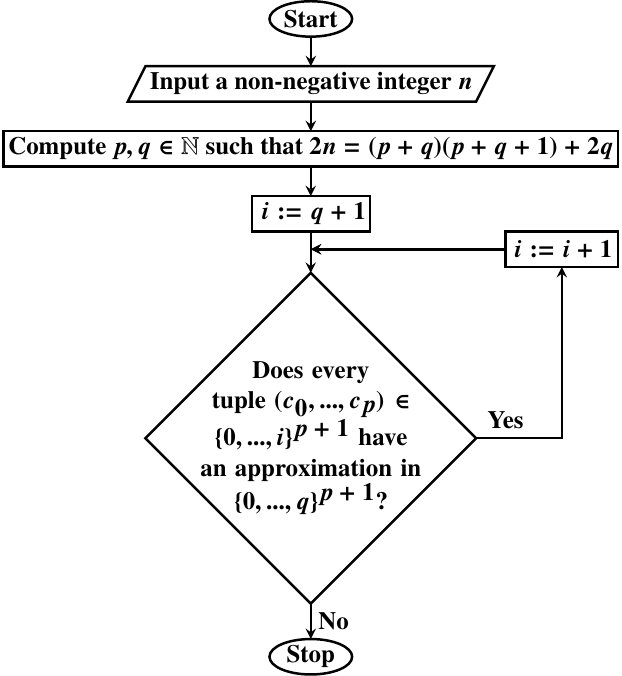}
\end{center}
\vskip 0.01truecm
\centerline{{\bf Flowchart 3}}
\vskip 0.1truecm
\centerline{A \mbox{semi-algorithm} which takes as input \mbox{$n \in \N$} and terminates if and only if \mbox{$n \in \N \setminus T$}}
\end{proof}
\begin{theorem}\label{theo8}
The set $T$ is not recursively enumerable.
\end{theorem}
\begin{proof}
It follows from Theorems~\ref{theo6} and~\ref{theo7}.
\end{proof}
\section{\normalsize A limit-computable function \mbox{$h:\N \to \N$} which eventually dominates
every computable function \mbox{$g:\N \to \N$}}
\begin{lemma}\label{lemm3} (\cite[p.~110]{Tyszka0}).
For non-negative integers, the equation \mbox{$x+y=z$} is equivalent to a system which consists of equations of the
forms \mbox{$v+1=w$} and \mbox{$u \cdot v=w$}.
\end{lemma}
\par
For \mbox{$n \in \N$}, \mbox{$h(n)$} denotes the smallest \mbox{$b \in \N$} such that
if a system of equations \mbox{$S \subseteq$}
\mbox{$\{x_j+1=x_k,~x_i \cdot x_j=x_k:~i,j,k \in \{0,\ldots,n\}\}$}
has a solution in \mbox{${\N}^{n+1}$},
then $S$ has a solution in \mbox{$\{0,\ldots,b\}^{n+1}$}.
From Lemma~\ref{lemm3} and \cite{Tyszka1}, it follows that
the function \mbox{$h:\N \to \N$} is computable in the limit and
eventually dominates every computable function \mbox{$g:\N \to \N$}.
\newpage
\section{\normalsize Two undecidable decision problems which are related to the function $h$}
\begin{theorem}\label{theo9}
No algorithm takes as input \mbox{non-negative integers} $n$ and $m$ and decides whether or not
\[
\forall (x_0,\ldots,x_n) \in {\N}^{n+1} ~\exists (y_0,\ldots,y_n) \in \{0,\ldots,m\}^{n+1}
\]
\[
((\forall j,k \in \{0,\ldots,n\} ~(x_j+1=x_k \Rightarrow y_j+1=y_k))~\wedge
\]
\[
(\forall i,j,k \in \{0,\ldots,n\} ~(x_i \cdot x_j=x_k \Rightarrow y_i \cdot y_j=y_k)))
\]
\end{theorem}
\begin{proof}
It holds because the function $h$ is not computable, and for every \mbox{$n \in \N$},
$h(n)$ is the smallest \mbox{$b \in \N$} such that
\[
\forall (x_0,\ldots,x_n) \in {\N}^{n+1} ~\exists (y_0,\ldots,y_n) \in \{0,\ldots,b\}^{n+1}
\]
\[
((\forall j,k \in \{0,\ldots,n\} ~(x_j+1=x_k \Rightarrow y_j+1=y_k))~\wedge
\]
\[
(\forall i,j,k \in \{0,\ldots,n\} ~(x_i \cdot x_j=x_k \Rightarrow y_i \cdot y_j=y_k)))
\]
\end{proof}
\begin{theorem}\label{theo10}
No algorithm takes as input a non-negative integer $n$ and decides whether or not
\[
\exists p,q \in \N ~((2n=(p+q)(p+q+1)+2q) \wedge
\]
\[
\forall (x_0,\ldots,x_p) \in {\N}^{p+1} ~\exists (y_0,\ldots,y_p) \in \{0,\ldots,q\}^{p+1}
\]
\[
((\forall j,k \in \{0,\ldots,p\} ~(x_j+1=x_k \Rightarrow y_j+1=y_k))~\wedge
\]
\[
(\forall i,j,k \in \{0,\ldots,p\} ~(x_i \cdot x_j=x_k \Rightarrow y_i \cdot y_j=y_k))))
\]
\end{theorem}
\begin{proof}
It follows from Theorem~\ref{theo9} and Lemma~\ref{Cantor}.
\end{proof}
\section{\normalsize A simple example of a non-recursively enumerable set \mbox{$W \subseteq \N$}
which is related to the function $h$}
Let $W$ denote the algorithmically undecidable subset of $\N$ considered in Theorem~\ref{theo10}.
Similarly as in Theorem~\ref{theo7}, the set \mbox{$\N \setminus W$} is recursively enumerable.
Similarly as in Theorem~\ref{theo8}, the set $W$ is not recursively enumerable.
Let \mbox{$\be:{\N}^3 \to \N$} denote G\"odel's $\be$ function, see \cite{Godel}.
For \mbox{$x_1,x_2,x_3 \in \N$},
\mbox{$\be(x_1,x_2,x_3)$} equals the remainder after integer division of $x_1$ by \mbox{$1+(x_3+1) \cdot x_2$}.
\begin{lemma}\label{Godel} (\cite{Godel}).
If \mbox{$(d_0,\ldots,d_p) \in {\N}^{p+1}$}, then $\exists z_1,z_2 \in \N ~\forall l \in \{0,\ldots,p\} ~\be(z_1,z_2,l)=d_l$.
\end{lemma}
\par
By Lemma~\ref{Godel}, the set $W$ consists of all \mbox{$n \in \N$} such that
\[
\forall u,v \in \N ~\exists a,b,p,q \in \N ~((2n=(p+q)(p+q+1)+2q) ~\wedge ~\forall i,j,k \in \{0,\ldots,p\}
\]
\[
((\be(a,b,i) \leqslant q) ~\wedge ~(\be(u,v,j)+1=\be(u,v,k) \Rightarrow \be(a,b,j)+1=\be(a,b,k)) ~\wedge
\]
\[
(\be(u,v,i) \cdot \be(u,v,j)=\be(u,v,k) \Rightarrow \be(a,b,i) \cdot \be(a,b,j)=\be(a,b,k))))
\]
In Peano arithmetic, the above formula expresses that
\[
\forall~u~v~\exists~a~b~p~q~((2n=(p+q)(p+q+1)+2q)~\wedge
\]
\[
\forall~i~j~k~(((i \leqslant p) \wedge (j \leqslant p) \wedge (k \leqslant p))~\Rightarrow
~\exists~x_1~x_2~x_3~\wt{x_1}~\wt{x_2}~\wt{x_3}~y_1~y_2~y_3~\wt{y_1}~\wt{y_2}~\wt{y_3}
\]
\[
((x_1<1+(i+1) \cdot v) \wedge (x_2<1+(j+1) \cdot v) \wedge (x_3<1+(k+1) \cdot v)~\wedge
\]
\[
(y_1<1+(i+1) \cdot b) \wedge (y_2<1+(j+1) \cdot b) \wedge (y_3<1+(k+1) \cdot b)~\wedge
\]
\[
(u=\wt{x_1} \cdot (1+(i+1) \cdot v)+x_1) \wedge (u=\wt{x_2} \cdot (1+(j+1) \cdot v)+x_2) \wedge (u=\wt{x_3} \cdot (1+(k+1) \cdot v)+x_3)~\wedge
\]
\[
(a=\wt{y_1} \cdot (1+(i+1) \cdot b)+y_1) \wedge (a=\wt{y_2} \cdot (1+(j+1) \cdot b)+y_2) \wedge (a=\wt{y_3} \cdot (1+(k+1) \cdot b)+y_3)~\wedge
\]
\[
(y_1 \leqslant q) \wedge (x_2+1=x_3 \Rightarrow y_2+1=y_3) \wedge (x_1 \cdot x_2=x_3 \Rightarrow y_1 \cdot y_2=y_3))))
\]
\par
No known set \mbox{$X \subseteq \N$} has an intuitive definition in Peano arithmetic
and is not recursively enumerable.
\section{\normalsize A simple example of a non-recursively enumerable set \mbox{$Y \subseteq \N$}, where\\
$Y=\{\frac{(p+q-2)(p+q-1)}{2}+q-1: (p,q \in \N \setminus \{0\}) \wedge \Phi(p,q)\}$ and the formula\\
$(p,q \in \N \setminus \{0\}) \wedge \Phi(p,q)$ concerns $((\N \setminus \{0\})^p,\{1,\ldots,q\}^p)$}
\begin{lemma}\label{lemm4}
For every \mbox{$x,y \in \N \setminus \{0\}$}, $x^2+y^2+xy+xy=(x+y)^2$.
For every \mbox{$x \in \N \setminus \{0\}$}, \mbox{$x=1$} if and only if \mbox{$x^2=x$}.
\end{lemma}
\begin{lemma}\label{lemm5}
In \mbox{$\N \setminus \{0\}$}, squaring, addition and using only conjunctions existentially defines multiplication and $1$.
\end{lemma}
\begin{proof}
It follows from Lemma~\ref{lemm4}.
\end{proof}
\begin{lemma}\label{lemm6}
In \mbox{$\N \setminus \{0\}$}, every polynomial equation $P(x_1,\ldots,x_n)=0$ with integer coefficients
is equivalent to a system of equations of the forms \mbox{$v^2=w$} and \mbox{$u+v=w$}.
\end{lemma}
\begin{proof}
It follows from Lemma~\ref{lemm5}.
\end{proof}
\begin{theorem}\label{theo11}
The set
\[
\{\frac{(p+q-2)(p+q-1)}{2}+q-1: (p,q \in \N \setminus \{0\}) ~\wedge
\]
\[
(\forall x_1,\ldots,x_p \in \N \setminus \{0\} ~\exists y_1,\ldots,y_p \in \{1,\ldots,q\} ~\forall i,j,k \in \{1,\ldots,p\}
\]
\[
((x_j^2=x_k \Rightarrow y_j^2=y_k) \wedge (x_i+x_j=x_k \Rightarrow y_i+y_j=y_k)))\}
\]
is co-recursively enumerable and not recursively enumerable.
\end{theorem}
\begin{proof}
It follows from Lemmas~\ref{Cantor} and \ref{lemm6}.
\end{proof}
\section{\normalsize A non-recursively enumerable set \mbox{$E \subseteq \N$}
whose definition uses exponentiation in \mbox{$\N \setminus \{0,1\}$}}
\begin{lemma}\label{lemm7}
For every \mbox{$x \in \N \setminus \{0,1\}$}, \mbox{$x=2$} if and only if
\[
\exists y,t \in \N \setminus \{0,1\} ~((x^{\textstyle y}=t) \wedge (y^{\textstyle x}=t) \wedge (x^{\textstyle x}=y))
\]
\end{lemma}
\begin{proof}
The above conjunction holds when \mbox{$(x=2) \wedge (y=4) \wedge (t=16)$}.
The above conjunction implies that
\[
x^{\textstyle x^{\textstyle x}}=x^{\textstyle y}=t=y^{\textstyle x}=\left(x^{\textstyle x}\right)^{\textstyle x}=x^{\textstyle x^2}
\]
Since $x \in \N \setminus \{0,1\}$, $x^{\textstyle x}=x^2$. Hence, $x=2$.
\end{proof}
\begin{lemma}\label{lemm8}
For every \mbox{$x,y,z \in \N \setminus \{0,1\}$},
\mbox{$x \cdot y=z$} if and only if \mbox{$(2^{\textstyle x})^{\textstyle y}=2^{\textstyle z}$}.
\end{lemma}
\begin{lemma}\label{lemm9}
For every \mbox{$x,y,z \in \N \setminus \{0,1\}$}, \mbox{$x+y=z$} if and only if
\mbox{$(2^{\textstyle 2^{\textstyle x}})^{\textstyle 2^{\textstyle y}}=2^{\textstyle 2^{\textstyle z}}$}.
\end{lemma}
\begin{lemma}\label{lemm10}
In \mbox{$\N \setminus \{0,1\}$}, exponentiation and using only conjunctions
existentially defines $2$, multiplication and addition.
\end{lemma}
\begin{proof}
It follows from Lemmas~\ref{lemm7}--\ref{lemm9}.
\end{proof}
\begin{lemma}\label{lemm11}
Every polynomial equation $P(x_1,\ldots,x_n)=0$ with integer coefficients
is equivalent to a system of equations of the forms \mbox{$2=w$}, \mbox{$u \cdot v=w$} and \mbox{$u+v=w$}.
\end{lemma}
\begin{proof}
We express the equation \mbox{$2 \cdot P(x_1,\ldots,x_n)=0$} as
\mbox{$P_1(x_1,\ldots,x_n)=P_2(x_1,\ldots,x_n)$}, where \mbox{$P_1(x_1,\ldots,x_n)$} and \mbox{$P_2(x_1,\ldots,x_n)$}
are polynomials whose coefficients are positive and even.
\end{proof}
\begin{lemma}\label{lemm12}
In \mbox{$\N \setminus \{0,1\}$}, every polynomial equation $P(x_1,\ldots,x_n)=0$
with integer coefficients is equivalent
to a system of equations of the form \mbox{$u^{\textstyle v}=w$}.
\end{lemma}
\begin{proof}
It follows from Lemmas~\ref{lemm10} and~\ref{lemm11}.
\end{proof}
\begin{theorem}\label{theo12}
The set
\[
\{2^{\textstyle p} \cdot 3^{\textstyle q}: (p,q \in \N \setminus \{0,1\}) ~\wedge
\]
\[
(\forall x_2,\ldots,x_p \in \N \setminus \{0,1\} ~\exists y_2,\ldots,y_p \in \{2,\ldots,q\}
~\forall i,j,k \in \{2,\ldots,p\} ~(x_i^{\textstyle x_j}=x_k \Rightarrow y_i^{\textstyle y_j}=y_k))\}
\]
is co-recursively enumerable and not recursively enumerable.
\end{theorem}
\begin{proof}
It follows from Lemma~\ref{lemm12}.
\end{proof}

\vskip 0.01truecm
\noindent
Apoloniusz Tyszka\\
Hugo Ko{\l}{\l}\k{a}taj University\\
Balicka 116B, 30-149 Krak\'ow, Poland\\
E-mail: \url{rttyszka@cyf-kr.edu.pl}
\end{sloppypar}
\end{document}